\documentclass[12pt]{article}

\begin{document}

\title{$sk$-Spline interpolation on ${\bf R}^{n}$ }
\author{F. Jarad, A. Kushpel, J. Levesley, K. Ta\c{s}}
\maketitle

\begin{center}
Abstract
\end{center}

The main aim of this article is to introduce $sk$-splines on ${\bf R}^{n}$ and establish representations of cardinal $sk$-splines with knots and points
of interpolation on the sets ${\bf A}{\bf Z}^{n}$, where ${\bf A}$
is an arbitrary nonsingular $n\times n$ matrix. Such sets of points are
analogs for ${\bf R}^{n}$ of number theoretic Korobov's grids on the
torus and proved to be useful for problems of very high dimensionality.

Keywords: $sk$-spline, interpolation.

Subject: 41A05, 42B05.

\section{Introduction}

$sk$-Splines were introduced and their basic theory was developed by Kushpel [2-4,7,8]. These results were
developed in [1,5,6,9]. In this cycle of works, $sk
$-splines were proved to be useful in Kolmogorov's problem on $n$-widths and
approximation of smooth, infinitely smooth, analytic and entire functions on
the circle, $d$-torus and $d$-dimensional sphere. The theory of $sk$-splines
and methods developed were applied by different authors for calculation of
sharp values of $n$-widths (see e.g. [11,13] and
references therein). Here we introduce $sk$-splines on ${\bf R}^{n}$. Let
$\Lambda :=\left\{ {\bf x}_{{\bf k}}\right\} $ be an additive group of
lattice points in ${\bf R}^{n}$ \ and $K\left( \cdot \right) :{\bf R}%
^{n}\rightarrow {\bf R}$ be a fixed kernel function. An $sk$-Spline is a
function of the form%
\[
sk\left( {\bf x}\right) =\sum_{{\bf x}_{{\bf k}}\in {\bf \Lambda
}}c_{{\bf k}}K\left( {\bf x-x}_{{\bf k}}\right) \in L_{2}({\bf R}^{n}).
\]%
Let us consider in more details the problem of interpolation in ${\bf R}%
^{n}$. For a fixed $\Lambda $ and any continuous function $f$ $:{\bf R}%
^{n}\rightarrow {\bf R}$ we need to find such $c_{{\bf j}}$, ${\bf j%
}\in {\bf Z}^{n}$ that
\[
\sum_{{\bf x}_{{\bf j}}\in \Lambda }c_{{\bf j}}K\left( {\bf x}_{%
{\bf k}}-{\bf x}_{{\bf j}}\right) =f\left( {\bf x}_{{\bf k}%
}\right) ,
\]%
for any ${\bf x}_{{\bf k}}\in \Lambda $. Assuming some general
conditions on the kernel function $K$ we \ give an explicit solution of this
problem in the case $\Lambda ={\bf A}{\bf Z}^{n}\subset {\bf R}^{n}$%
, where ${\bf A}$ is an arbitrary nonsingular matrix. Such grids are
analogs for ${\bf R}^{n}$ of number theoretic Korobov's grids and their
various modifications, such as sparse grids. These grids proved to be useful
in high dimensional integration and interpolation (see e.g. [1]).

Let $L_{p}({\bf R}^{n})$ be the usual space of $p$-integrable functions
equipped with the norm
\[
\Vert f\Vert _{p}=\Vert f\Vert _{L_{p}({\bf R}^{n})}:=\left\{
\begin{array}{cc}
\left( \int_{{\bf R}^{n}}\left\vert f({\bf x})\right\vert ^{p}d{\bf %
x}\right) ^{1/p}, & 1\leq p<\infty , \\
\mathrm{ess}\,\,\sup_{{\bf x}\in {\bf R}^{n}}|f({\bf x})|, &
p=\infty .%
\end{array}%
\right.
\]
Let ${\bf x}$ and ${\bf y}$ be two vectors in ${\bf R}^{n}$ and $%
\left\langle {\bf x,y}\right\rangle =\sum_{k=1}^{n}{\bf x}_{k}{\bf y}_{k}$
be the usual scalar product, $|{\bf x}|=\left\langle {\bf x,x}\right\rangle^{1/2}$. For an integrable on ${\bf R}^{n}$ function,
i.e., $f({\bf x})\in L_{1}\left( {\bf R}^{n}\right) $ define its
Fourier transform
\[
{\bf F}f({\bf y})=\int_{{\bf R}^{n}}\exp \left( -i\left\langle
{\bf x,y}\right\rangle \right) f({\bf x})d{\bf x}.
\]%
and its formal inverse as
\[
\left( {\bf F}^{-1}f\right) ({\bf x})=\frac{1}{\left( 2\pi \right) ^{n}%
}\int_{{\bf R}^{n}}\exp \left( i\left\langle {\bf x,y}\right\rangle
\right) f({\bf y})d{\bf y}.
\]
We will need the following well-known results (see [12] for more
information).

{\bf Theorem 1} \begin{em}
(Plancherel's theorem) The Fourier transform is a
linear continuous operator from $L_{2}\left( {\bf R}^{n}\right) $ onto $%
L_{2}\left( {\bf R}^{n}\right) .$ The inverse Fourier transform, $F^{-1},$
can be obtained by letting
\[
\left( {\bf F}^{-1}g\right) \left( {\bf x}\right) =\frac{1}{(2\pi )^{n}%
}\left( {\bf F}g\right) \left( -{\bf x}\right)
\]%
for any $g\in L_{2}\left( {\bf R}^{n}\right) .$ 
\end{em}

{\bf Theorem 2} 
\begin{em}
(Poisson summation formula) Suppose that for some $C>0$ and $\delta >0$ we have $\max \left\{ f({\bf x}),{\bf F}f(%
{\bf x})\right\} \leq C\left( 1+\left\vert {\bf x}\right\vert \right)
^{-n-\delta }$. Then
\[
\sum_{{\bf m}\in {\bf Z}^{n}}f\left( {\bf x}+{\bf m}\right)
=\sum_{{\bf m}\in {\bf Z}^{n}}{\bf F}f\left( 2\pi {\bf m}\right)
\exp \left( 2\pi i\left\langle {\bf m},{\bf x}\right\rangle \right)
\]%
and the series converges absolutely. 
\end{em}

\section{Interpolation by $sk$-splines on ${\bf R}^{n}$}

Let ${\bf A}$ be a nonsingular $n\times n$ matrix. Consider the set $%
\Omega _{{\bf A}}$ of grid points ${\bf x}_{{\bf m}}:={\bf Am}$,
where ${\bf m }\in{\bf Z}^{n}$. For a fixed continuous kernel function
$K$, the space $SK\left( \Omega _{{\bf A}}\right) $ of $sk$-splines on $%
\Omega _{{\bf A}}$ is the space of functions representable in the form
\[
sk\left( {\bf x}\right) =\sum_{{\bf m }\in{\bf Z}^{n}}c_{{\bf m}%
}K\left( {\bf x-x}_{{\bf m}}\right) ,
\]%
where $c_{{\bf m}}\in {\bf R}$, ${\bf m\in }{\bf Z}^{n}$. Let $f\left( {\bf x}%
\right) $ be a continuous function, $f:{\bf R}^{n}\rightarrow {\bf R}$%
. Consider the problem of interpolation by $sk$-splines, $sk\left( {\bf x}%
_{\bf s}\right) =f\left( {\bf x}_{\bf s}\right) $, where $%
{\bf s }\in{\bf Z}^{n}$. Even in the one-dimensional case the problem
of interpolation does not always have a solution. If the solution exists then
the $sk$-spline interpolant can be written in the form
\[
sk\left( {\bf x}\right) =\sum_{{\bf s }\in{\bf Z}^{n}}f\left(
{\bf x}_{{\bf s}}\right) \widetilde{sk}\left( {\bf x}-{\bf x}_{%
{\bf s}}\right) ,
\]%
where $\widetilde{sk}\left( {\bf x}\right) $ is a fundamental $sk$%
-spline, i.e.%
\[
\widetilde{sk}\left( {\bf x}_{{\bf s}}\right) =\left\{
\begin{array}{cc}
1, & {\bf s=0,} \\
0, & {\bf s\neq 0.}%
\end{array}%
\right.
\]

{\bf Theorem 3} 
\begin{em}
Let ${\bf A}$ be a
nonsingular matrix, $K:{\bf R}^{n}\rightarrow {\bf R}$,$K\in
L_{2}\left( {\bf R}^{n}\right) \cap L_{1}\left( {\bf R}^{n}\right)
\cap C\left( {\bf R}^{n}\right) .$ Assume that
\[
\sum_{{\bf m}\in {\bf Z}^{n}}{\bf F}\left( K\right) \left( {\bf%
z+}2\pi \left( {\bf A}^{-1}\right) ^{T}{\bf m}\right) \neq 0,\forall
{\bf z }\in2\pi \left( {\bf A}^{-1}\right) ^{T}{\bf Q},
\]%
where ${\bf Q}:=\left\{ {\bf x|x=}\left( x_{1},\cdot \cdot \cdot
,x_{n}\right) \in {\bf R}^{n},0\leq x_{k}\leq 1,1\leq k\leq
n\right\} ,$ and the function
\[
\Upsilon \left( {\bf z}\right) :=\frac{1}{\sum_{{\bf m}\in {\bf Z}%
^{n}}{\bf F}\left( K\right) \left( {\bf z+}2\pi \left( {\bf A}%
^{-1}\right) ^{T}{\bf m}\right) }
\]%
can be represented by its Fourier series, i.e. for any ${\bf z\in }%
{\bf R}^{n}$, \
\begin{equation}
\Upsilon \left( {\bf z}\right) =\sum_{{\bf s}\in {\bf Z}^{n}}\alpha
_{{\bf s}}\exp \left( -i\left\langle {\bf As},{\bf z}\right\rangle
\right) .  \label{repr}
\end{equation}%
Then%
\[
\widetilde{sk}\left( {\bf x}\right) =\frac{\det \left( {\bf A}\right)
}{\left( 2\pi \right) ^{n}}\int_{{\bf R}^{n}}\Upsilon \left( {\bf z}%
\right) {\bf F}\left( K\right) \left( {\bf z}\right) \exp \left(
i\left\langle {\bf z},{\bf x}\right\rangle \right) d{\bf z}
\]%
and this representation is unique. 
\end{em}

{\bf Proof} Since $K\in L_{2}\left( {\bf R}^{n}\right) \cap
L_{1}\left( {\bf R}^{n}\right) $ and \ using (\ref{repr}) we get
\[
\widetilde{sk}\left( {\bf x}\right) =\frac{\det \left( {\bf A}\right)
}{\left( 2\pi \right) ^{n}}\int_{{\bf R}^{n}}\Upsilon \left( {\bf z}%
\right) {\bf F}\left( K\right) \left( {\bf z}\right) \exp \left(
i\left\langle {\bf x},{\bf z}\right\rangle \right) d{\bf z}
\]%
\[
=\frac{\det \left( {\bf A}\right) }{\left( 2\pi \right) ^{n}}\int_{%
{\bf R}^{n}}{\bf F}\left( K\right) \left( {\bf z}\right) \left(
\sum_{{\bf s}\in {\bf Z}^{n}}\alpha _{{\bf s}}\exp \left(
-i\left\langle {\bf As},{\bf z}\right\rangle \right) \right) \exp
\left( i\left\langle {\bf x,z}\right\rangle \right) d{\bf z}
\]%
\[
=\frac{\det \left( {\bf A}\right) }{\left( 2\pi \right) ^{n}}\sum_{%
{\bf s}\in {\bf Z}^{n}}\alpha _{\bf s}\int_{{\bf R}^{n}}%
{\bf F}\left( K\right) \left( {\bf z}\right) \exp \left( i\left\langle
{\bf x}-{\bf As},{\bf z}\right\rangle \right) d{\bf z.}
\]%
Since $K\in L_{2}\left( {\bf R}^{n}\right) $ then by Theorem 1%
\[
\widetilde{sk}\left( {\bf x}\right) =\det \left( {\bf A}\right) \sum_{%
{\bf s}\in {\bf Z}^{n}}\alpha _{{\bf s}}K\left( {\bf x}-{\bf
As}\right) ,
\]%
so that $\widetilde{sk}\left( {\bf x}\right) \in SK\left( \Omega _{%
{\bf A}}\right) $. Next, we calculate $\widetilde{sk}\left( {\bf Am}%
\right) $ for ${\bf m}\in {\bf Z}^{n}$,
\[
\widetilde{sk}\left( {\bf Am}\right) =\frac{\det \left( {\bf A}\right)
}{\left( 2\pi \right) ^{n}}\int_{{\bf R}^{n}}\Upsilon \left( {\bf z}%
\right) {\bf F}\left( K\right) \left( {\bf z}\right) \exp \left(
i\left\langle {\bf z},{\bf Am}\right\rangle \right) d{\bf z}
\]%
\[
=\frac{\det \left( {\bf A}\right) }{\left( 2\pi \right) ^{n}}\sum_{%
{\bf l}\in {\bf Z}^{n}}\int_{2\pi \left( {\bf A}^{-1}\right)^{T}%
{\bf l}+2\pi \left( {\bf A}^{-1}\right) ^{T}{\bf Q}}\Upsilon \left(
{\bf z}\right) {\bf F}\left( K\right) \left( {\bf z}\right) \exp
\left( i\left\langle {\bf z},{\bf Am}\right\rangle \right) d{\bf z}
\]%
\[
=\frac{\det \left( {\bf A}\right) }{\left( 2\pi \right) ^{n}}
\]%
\[
\times \sum_{{\bf l}\in {\bf Z}^{n}}\int_{2\pi \left( {\bf A}%
^{-1}\right) ^{T}{\bf Q}}\frac{{\bf F}\left( K\right) \left( {\bf z-%
}2\pi \left( {\bf A}^{-1}\right) ^{T}{\bf l}\right)
}{\sum_{{\bf m}\in {\bf Z}^{n}}{\bf F}\left( K\right) \left(
{\bf z-}2\pi \left( {\bf A}^{-1}\right) ^{T}{\bf l}%
+2\pi \left( {\bf A}^{-1}\right) ^{T}{\bf m}\right) }
\]%
\[
\times \exp \left( i\left\langle {\bf z-}2\pi \left( {\bf A}%
^{-1}\right) ^{T}{\bf l},{\bf Am}\right\rangle \right) d%
{\bf z.}
\]%
Changing ${\bf l}$ by $-{\bf l}$ in $\sum_{{\bf l}\in {\bf Z}%
^{n}}$ we get%
\[
\widetilde{sk}\left( {\bf Am}\right) =\frac{\det \left( {\bf A}\right)
}{\left( 2\pi \right) ^{n}}
\]%
\[
\times \sum_{{\bf l}\in {\bf Z}^{n}}\int_{2\pi \left( {\bf A}%
^{-1}\right) ^{T}{\bf Q}}\frac{{\bf F}\left( K\right) \left( {\bf z+%
}2\pi \left( {\bf A}^{-1}\right) ^{T}{\bf l}\right)
}{\sum_{{\bf m}\in {\bf Z}^{n}}{\bf F}\left( K\right) \left(
{\bf z+}2\pi \left( {\bf A}^{-1}\right) ^{T}{\bf l}%
+2\pi \left( {\bf A}^{-1}\right) ^{T}{\bf m}\right) }
\]%
\[
\times \exp \left( i\left\langle {\bf z+}2\pi \left( {\bf A}%
^{-1}\right) ^{T}{\bf l},{\bf Am}\right\rangle \right) d%
{\bf z.}
\]%
Since
\[
\sum_{{\bf m}\in {\bf Z}^{n}}{\bf F}\left( K\right) \left( {\bf%
z+}2\pi \left( {\bf A}^{-1}\right) ^{T}{\bf l}+2\pi
\left( {\bf A}^{-1}\right) ^{T}{\bf m}\right)
\]%
\[
=\sum_{{\bf m}\in {\bf Z}^{n}}{\bf F}\left( K\right) \left( {\bf%
z+}2\pi \left( {\bf A}^{-1}\right) ^{T}{\bf m}\right)
=\Upsilon ^{-1}\left( {\bf z}\right)
\]%
and
\[
\exp \left( i\left\langle {\bf z+}2\pi \left( {\bf A}%
^{-1}\right) ^{T}{\bf l},{\bf Am}\right\rangle \right) =\exp
\left( i\left\langle {\bf z},{\bf Am}\right\rangle \right)
\]%
for any ${\bf l}\in {\bf Z}^{n}$ then
\[
\widetilde{sk}\left( {\bf Am}\right) =\frac{\det \left( {\bf A}\right)
}{\left( 2\pi \right) ^{n}}\sum_{{\bf l}\in {\bf Z}^{n}}\int_{2\pi
\left( {\bf A}^{-1}\right) ^{T}{\bf Q}}\Upsilon \left( {\bf z}%
\right) {\bf F}\left( K\right) \left( {\bf z+}2\pi \left( {\bf A}%
^{-1}\right) ^{T}{\bf l}\right)
\]%
\[
\times \exp \left( i\left\langle {\bf z+}2\pi \left( {\bf A}%
^{-1}\right) ^{T}{\bf l},{\bf Am}\right\rangle \right) d%
{\bf z}
\]%
\[
=\frac{\det \left( {\bf A}\right) }{\left( 2\pi \right) ^{n}}\int_{2\pi
\left( {\bf A}^{-1}\right) ^{T}{\bf Q}}\exp \left( i\left\langle
{\bf z},{\bf Am}\right\rangle \right) \Upsilon \left( {\bf z}%
\right) \sum_{{\bf l}\in {\bf Z}^{n}}{\bf F}\left( K\right) \left(
{\bf z}+2\pi \left( {\bf A}^{-1}\right) ^{T}{\bf l}\right)
d{\bf z}
\]%
\[
=\frac{\det \left( {\bf A}\right) }{\left( 2\pi \right) ^{n}}\int_{2\pi
\left( {\bf A}^{-1}\right) ^{T}{\bf Q}}\exp \left( i\left\langle
{\bf z},{\bf Am}\right\rangle \right) d{\bf z}
\]%
\[
=\int_{{\bf Q}}\exp \left( i\left\langle {\bf z},{\bf m}%
\right\rangle \right) d{\bf z}=\left\{
\begin{array}{cc}
1, & m_{k}=0,1\leq k\leq n, \\
0, & otherwise.%
\end{array}%
\right. .
\]
The representation of the fundamental $sk-$spline is unique since the
functions $\widetilde{sk}\left( {\bf x}-{\bf x}_{{\bf m}}\right) ,%
{\bf x}_{{\bf m}}\in \Omega _{{\bf A}}$ are linearly independent.  

Observe that Theorem 3 is a multidimensional analog of Schoenberg's result
on interpolation by polynomial splines of odd degree on ${\bf R}$ [10].
%
%

Let ${\bf Q}$ be a nonsingular matrix and $f\in L_{2}\left(
{\bf R}^{n}\right) $. Then%
\[
{\bf F}\left( f\left( {\bf Q}\cdot \right) \right) \left( {\bf x}%
\right) =\frac{1}{\det \left( {\bf Q}\right) }{\bf F}\left( f\left(
\cdot \right) \right) \left( \left( {\bf Q}^{-1}\right) ^{T}{\bf x}%
\right) .
\]%
Hence for any $K\in L_{2}\left( {\bf R}^{n}\right) $, by Theorem 1,
\[
{\bf F}\left( {\bf F}\left( K\right) \left( 2\pi \left( {\bf A}%
^{-1}\right) ^{T}\left(\cdot \right) \right) \right) \left(
{\bf y}\right)
\]%
\[
=\frac{1}{\det \left( 2\pi \left( {\bf A}^{-1}\right) ^{T}\right) }%
{\bf F\circ F}\left( K\left( \cdot \right) \right) \left( \left( 2\pi
\left( {\bf A}^{-1}\right) ^{T}\right) ^{-1}{\bf y}\right)
\]%
\[
=\frac{\det \left( {\bf A}\right) }{\left( 2\pi \right) ^{n}}\left( 2\pi
\right) ^{n}{\bf F}^{-1}\circ { \bf F}\left( K\left( \cdot \right)
\right) \left( -\left( 2\pi \left( {\bf A}^{-1}\right) ^{T}\right) ^{-1}%
{\bf y}\right)
\]%
\begin{equation}
=\det \left( {\bf A}\right) K\left( -\frac{1}{2\pi }{\bf A}^{T}{\bf
y}\right) .  \label{sk111}
\end{equation}%
Assume that $K$ satisfies the conditions of Theorem 2. Let ${\bf z}=2\pi
\left( {\bf A}^{-1}\right) ^{T}{\bf x}$ then using (\ref{sk111}) we get%
\[
\Upsilon ^{-1}\left( {\bf z}\right) =\sum_{{\bf m}\in {\bf Z}^{n}}%
{\bf F}\left( K\right) \left( {\bf z+}2\pi \left( {\bf A}%
^{-1}\right) ^{T}{\bf m}\right)
\]%
\[
=\sum_{{\bf m}\in {\bf Z}^{n}}{\bf F}\left( K\right) \left( 2\pi
\left( {\bf A}^{-1}\right) ^{T}\left( {\bf x}+{\bf m}\right)
\right)
\]%
\[
=\det \left( {\bf A}\right) \sum_{{\bf m}\in {\bf Z}^{n}}K\left( -%
{\bf A}^{T}{\bf m}\right) \exp \left( 2\pi i\left\langle {\bf m,x}%
\right\rangle \right)
\]%
\[
=\det \left( {\bf A}\right) \sum_{{\bf m}\in {\bf Z}^{n}}K\left( -%
{\bf A}^{T}{\bf m}\right) \exp \left( 2\pi i\left\langle {\bf m,}%
\left( 2\pi \left( {\bf A}^{-1}\right) ^{T}\right) ^{-1}{\bf z}%
\right\rangle \right)
\]%
\[
=\det \left( {\bf A}\right) \sum_{{\bf m}\in {\bf Z}^{n}}K\left( -%
{\bf A}^{T}{\bf m}\right) \exp \left( i\left\langle {\bf Am,z}%
\right\rangle \right) .
\]%
Hence we proved

{\bf Corollary 1} 
\begin{em}
Let $K$ satisfies the conditions of
Theorem 1-Theorem 3 then
\[
\widetilde{sk}\left( {\bf x}\right) ={\bf F}^{-1}\left( \frac{{\bf F%
}\left( K\right) \left( \cdot \right) }{\sum_{{\bf m}\in {\bf Z}%
^{n}}K\left( -{\bf A}^{T}{\bf m}\right) \exp \left( i\left\langle
{\bf Am,\cdot }\right\rangle \right) }\right) \left( {\bf x}\right) .
\]
\end{em}

{\bf Example 1} Let ${\bf B}$ be a nonsingular matrix and $K\left(
{\bf x}\right) $ be a Gaussian of the form $K\left( {\bf x}\right)
=\exp \left( -\left\vert {\bf Bx}\right\vert ^{2}\right) ,$ then
\[
{\bf F}\left( K\right) \left( {\bf z}\right) =\frac{\pi ^{n/2}}{\det
\left( {\bf B}\right) }\exp \left( -\left\vert \left( {\bf B}%
^{-1}\right) ^{T}{\bf z}\right\vert ^{2}\right) .
\]%
It is easy to check that in this case $K$ satisfies the conditions of
Theorem 1-Theorem 3. Hence for the interpolation by ${\bf A}{\bf Z}^{n}
$-shifts of Gaussians we have
\[
\widetilde{sk}\left( {\bf x}\right) =\frac{\pi ^{n/2}}{\det \left(
{\bf B}\right) }{\bf F}^{-1}\left( \frac{\exp \left( -\left\vert
\left( {\bf B}^{-1}\right) ^{T}\cdot \right\vert ^{2}\right) }{%
\sum_{{\bf m}\in {\bf Z}^{n}}\exp \left( -\left\vert {\bf BA}^{T}%
{\bf m}\right\vert ^{2}+i\left\langle {\bf Am,\cdot }\right\rangle
\right) }\right) \left( {\bf x}\right) .
\]

\end{document}